
\input amstex
\documentstyle{amsppt}
\magnification1100

\def\s{\Cal S}
\def\L{\text{\bf L}}
\def\Z{\text{\bf Z}}

\topmatter
\title Controlled surgery with trivial local fundamental groups
  \endtitle

\author Erik Pedersen, Frank Quinn and Andrew Ranicki\endauthor
\rightheadtext{Controlled Surgery with trivial local fundamental groups}
\leftheadtext{Pedersen, Quinn and Ranicki}
\address
Dept. of Mathematical Sciences\newline
\indent Binghamton University\newline
\indent Binghamton, New York 13901, USA
\endaddress
\email erik\@math.binghamton.edu\endemail
\address 
Dept. of Mathematics\newline
\indent Virginia Tech\newline
\indent Blacksburg VA 24061-0123, USA
\endaddress
\email quinn\@math.vt.edu\endemail
\address 
Dept. of Mathematics and Statistics\newline
\indent University of Edinburgh\newline
\indent Edinburgh EH9 3JZ, Scotland, UK
\endaddress

\email aar\@maths.ed.ac.uk\endemail



\thanks Partially supported by the US National Science Foundation, 
EPSRC Grant GR/M82158/01 and the Leverhulme Trust. \endthanks

\abstract
We provide a proof of the controlled surgery sequence, including
stability, in the special case that the local fundamental groups are
trivial. Stability is a key ingredient in the construction of exotic
homology manifolds by Bryant, Ferry, Mio and Weinberger, but no proof has
been available. The development given here is based on work of 
M. Yamasaki.
\endabstract

\endtopmatter
In this note we provide a proof of the controlled surgery exact
sequence used in the construction of exotic homology manifolds by
Bryant, Ferry, Mio and Weinberger, \cite{BFMW}.  A primitive version of
controlled surgery was developed by the second author in his definition
of the invariant that identifies exotic homology manifolds, \cite{Q3,Q4}.  
Surgery with bounded control, including exact sequences, was
developed in \cite{FP}.  The remarkable limit construction of
\cite{BFMW} uses a refinement of the sequences of \cite{FP}.  Roughly
speaking \cite{FP} describes a limit as $\epsilon\to 0$ while
\cite{BFMW} depends crucially on a stability property of the limiting
process.  The proof of the refinement was postponed to a planned
project that was never completed.  The intent was to deduce stability
in general from a special case with an independent solution known to be
stable, the ``$\alpha$ approximation theorem'' of \cite{CF}.  This is
reasonable in principle and may be possible, but it has become clear
that the authors of \cite{BFMW} have not addressed serious technical
issues needed to actually carry it out.  As noted in the review
\cite{R3} {\it Until\/} [the planned project] {\it or some appropriate
substitute becomes available the surgery classification must be
regarded as somewhat provisional -- although there is little doubt
among the experts that it is correct\/}.  The purpose of this paper is
to provide the ``appropriate substitute''.  Our proof is direct, and is
based on work of Yamasaki \cite{Y}.

The following is Theorem 2.4 of \cite{BFMW} 
with minor inaccuracies corrected.

\proclaim{1. Theorem} 
Suppose $B$ is a finite dimensional compact metric ANR and a dimension
$n\geq4$ is given.  There is a stability threshold $\epsilon_0>0$ so
that for any $\epsilon_0>\epsilon>0$ there is $\delta>0$ with the
following property: If $f\:N\to B$ is $(\delta,1)$-connected and $N$ is
a compact $n$-manifold then there is a controlled surgery exact sequence
$$H_{n+1}(B,\L)@>>>\s_{\epsilon,\delta}(N,f)@>>>[N,\partial N;
G/TOP,*]@>>> H_n(B;\L).$$ 
\endproclaim

\subhead 2. Definitions\endsubhead
\roster
\item $\s_{\epsilon,\delta}(N,f)$ is the {\it controlled structure
set\/} and is the object of interest.  It is defined as the set of
equivalence classes of $(M,g)$, where $M$ is a topological manifold and
$g\:M\to N$ is a homeomorphism on boundaries and a $\delta$ homotopy
equivalence rel boundary.  The equivalence relation is given by:
$(M,g)\sim (M',g')$ if there is a homeomorphism $h\:M\to M'$ whose
restriction to the boundary commutes with $g$ and $g'$, and which
$\epsilon$ homotopy commutes on all of $M$.  It is part of the
assertion that this does actually give an equivalence relation.
\item Sizes refer to measurements in $B$.  For instance an $\epsilon$
homotopy between maps $M\to N$ is $G\:M\times I\to N$ so that for any
$x\in M$, the arc $fG(\{x\}\times I)$ lies in the ball of radius
$\epsilon$ about $fG(x,0)$.  For basic definitions see Quinn \cite{Q1}. 
\item $[N,\partial N; G/TOP]$ is the set of homotopy classes of maps of
$N$, rel boundary, to the classifying space $G/TOP$.  This is the same
term as appears in the uncontrolled surgery sequence.   
\item $H_*(B;\L)$ is homology with coefficients in the 4-periodic
simply-connected surgery spectrum $\L$, with
$\pi_*(\L)=L_*(\Z)$ for $*\geq 0$. The
spectrum can be constructed algebraically using (uncontrolled)
quadratic Poincar\'e chain complexes over the ring $\Z$ (\cite{R2}).
\endroster

\subhead 3. Outline\endsubhead
\medskip

The proof divides into two parts.  First there is a sequence with more
primitive obstruction terms:
$$L_{n+1}(B;\Z,\epsilon,\delta)@>>>
\s_{\epsilon,\delta}(N,f)@>>>[N,\partial N; G/TOP,*]@>>> 
L_n(B;\Z,\epsilon,\delta).$$
This is ``well known'' and quickly assembled from pieces in the
literature, though not entirely straightforward.  The second part shows
that an assembly map defined by Yamasaki,
$$H_n(B;\L)@>>>L_n(B;\Z,\epsilon,\delta)$$
is an isomorphism for suitable $\epsilon$, $\delta$.  Theorem 1 follows
by using this to replace the obstruction term in the primitive
sequence.  The proof of assembly isomorphism follows that of Yamasaki
for the $L^{-\infty}$ case, with the substitution of a $K$-theory
splitting argument for stabilization by $\times S^1$.

We remark that the key feature needed in constructing homology
manifolds is ``stability'': the sequence in Theorem 1 holds for
particular $\epsilon$ and $\delta$ rather than just for the
$\epsilon\to 0$ limit, even though the homotopy and obstruction terms
are explicitly independent of $\epsilon$ and $\delta$.  Stability in
the sequence comes from stability in the assembly isomorphism. 
Yamasaki's arguments are stable for relatively straightforward reasons,
so the basic source traces back through the $K$-theory splitting
argument to stability of vanishing of controlled Whitehead groups in
\cite{Q1, \S8}.  Finally we note that for stability the use of
$\epsilon$ control, rather than bounded or continuous control, 
is essential.  Bounded versions can be used to identify the
inverse limit as $\epsilon \to 0$, but do not give information on how
the limit is approached.

\subhead 4. Surgery\endsubhead

We describe how the ``primitive'' surgery sequence of \S3 is a
straightforward controlled version of the standard sequence.  The only
awkward point is that the surgery obstruction $[N,\partial N;G/TOP]\to
L_n(B;\Z,\epsilon,\delta)$ still must be defined using surgery rather
than by a direct chain-level construction.

  $L_n(B;\Z,\epsilon,\delta)$ is defined to be the group of
$n$-dimensional quadratic Poincar\'e complexes ($\Z$ coefficients,
\cite{R1}), over $B$ and with radius $<\delta$, modulo bordism through
$(n+1)$-dimensional quadratic Poincar\'e pairs with radius $<\epsilon$. 
Here we use the routine version of controlled algebra \cite{Q1} that
locates bases for modules at points in $B$, and measures radii of
ordinary homomorphisms in terms of distance between involved basis
points.  Yamasaki goes further in using ``geometric'' morphisms that
incorporate paths in a space.  This refinement is unnecessary here
because the local fundamental groups are constant (in fact trivial).

Suppose $M\to N$ is a normal map representing an element in
$[N,\partial N;G/TOP]$.  For any $\delta>0$ we can do surgery below the
middle dimension \cite{Q3} to make $M\to N$ $(\delta,j)$-connected over
$B$, where $n=2j$ or $2j+1$.  In this case we can give the relative
chains a $\delta$ controlled quadratic Poincar\'e structure:
nondegenerate forms with a symmetry condition in even dimensions, and
``short odd complexes'' \cite{R4,\S6} in odd dimensions.  We define the surgery
obstruction function by taking a normal map to the quadratic Poincar\'e
structure on a highly-connected normal map in the bordism class.

The verification that this function is well-defined uses a relative
construction.  If there is a normal bordism between highly-connected
normal maps then do surgery to make the bordism also $\delta$-connected
below the middle dimension.  The relative chains of the bordism can
then be given the structure of a $\delta$ quadratic Poincar\'e chain
bordism between the chains of the two maps.

The basis for exactness of the sequence is: if the quadratic structure
on the chains of a highly connected normal map is null-bordant through
a {\it highly-connected\/} quadratic chain bordism, then we get the
usual algebraic data for doing middle-dimensional surgery to get an
equivalence.  More specifically if $\epsilon>0$ then there is
$\delta>0$ so that if the quadratic chains of a highly
$\delta$-connected normal map is highly $\delta$-connected
algebraically nullbordant, then the normal map is bordant (by surgery)
to an $\epsilon$ equivalence.  The $\epsilon$ and $\delta$ here come
from the controlled Hurewicz and Whitehead theorems \cite{Q1, \S5}
rather than subtle stability issues.

This is the point at which dimension issues arise.  Standard surgery
requires dimension $\geq 5$, and gives surgery sequences for smooth and
$PL$ manifolds when the structure set and homotopy terms are changed
appropriately.  The topological version holds in dimension 4: since the
local fundamental groups are trivial \cite{Q2} gives the controlled
embeddings of 2-spheres in 4-manifolds needed for the surgery.  In fact
there is a weak 3-dimensional version in which objects in the structure
set are $\delta$ {\it homology\/} equivalences and equivalences are
$\epsilon$ homology $s$-cobordisms, but we have not tried to include this
in the main statement.

The algebra and topology are brought together by algebraic surgery on
controlled quadratic complexes.  The results needed are exactly
analogous to the topological case, and easily obtained by adding
control to the arguments of \cite{R1}: given $n$ and $\epsilon>0$ there
is $\delta$ so that an $n$-dimensional $\delta$ quadratic Poincar\'e
complex over $B$ is $\epsilon$ bordant to one that is connected up to
the middle dimension.  Similarly if two highly-connected quadratic
complexes are $\delta$-bordant then there is a highly connected
$\epsilon$-bordism.

\subhead 5. Assembly\endsubhead

The assembly map $H_n(B;\L)@>>>L_n(B;\Z,\epsilon,\delta)$ is
defined by Yamasaki \cite{Y, \S3} using ``cycles''.  The
characterization of assembly maps in Weiss and Williams \cite{WW} shows
this agrees with other definitions.  We review \cite{Y} to explain how
it reduces our case to a $K$-theory splitting problem.

An element of $H_n(B;\L)$ is represented by (1) a triangulated
codimension 0 submanifold $V$ of some sphere $S^N$; (2) a quadratic
Poincar\'e pair over each simplex of $V$ so that the pairs over
$\partial \sigma$ fit together to give the boundary of the pair over
$\sigma$, and so that they all glue together to give a Poincar\'e
complex of dimension $n$; and (3) a map $V\to B$.  The complex obtained
by glueing them all together is the assembly.  To get control on the
assembled complex one subdivides $V$ so finely that images of simplices
have diameter $<\delta$ in $B$.  Define a function from the basis of
the assembled complex into $B$ as follows: each basis element in the
assembled complex comes from one of the fragments, lying over some
simplex.  Take the basis element to an arbitrary point in the image of
the simplex.  Since the structure maps in the assembled complex keep
this basis element inside its fragment, therefore over points in a
small image, the assembled complex has radius $<\delta$.

This description of the assembly makes it clear that to show it is an
isomorphism we need to start with a controlled complex and split it up
as a union of small pieces.  Full splitting follows inductively from
being able to split once, so we are reduced to the following analog of
\cite{Y, Lemma 2.5}:

\proclaim{6.  Lemma: splitting Poincar\'e complexes} Suppose $B$ is a
metric space, $W\subset B$, $\epsilon>0$, and a dimension $n$ are
given.  Then there is $\delta>0$ so that an $n$-dimensional quadratic
Poincar\'e complex $D$ over $B$ with radius $<\delta$ is $\epsilon$
equivalent to a union of $\epsilon$ Poincar\'e pairs $D'\cup_CD''$ with
$D'$ located over $B-W$, $D''$ located over $W^{\epsilon}$, and the
common boundary $C$ located over $W^{\epsilon}\cap(B-W)$.\endproclaim

Yamasaki comes close to producing such a decomposition.  He produces
pairs $C\to D'$, $C\to D''$ with $D'$, $D''$ over $B-W$, $W^{\epsilon}$
respectively, as desired, but $C$ may be nontrivial over all of $B$. 
In a nutshell, $D''$ is the maximal based subcomplex of $D$ lying over
$B-W$, $D'$ is the quotient $D/D''$ with basis the complement of the
basis of $D''$ in the basis of $D$.  Since these are $n$-dimensional
complexes of radius $< \delta$, $D'$ lies over $W^{n\delta}$.  Finally
$C$ is Ranicki's ``algebraic boundary'' defined as the mapping cone of
a duality homomorphism.  This homomorphism uses all of $D$ so has no
restrictions on its location.  However the duality homomorphism is a
$\delta$ chain equivalence outside $W^{n\delta}\cap(W^{-n\delta})$, so
the mapping cone $C$ is $\delta$ contractible there.  Over the region
where it is contractible Yamasaki uses the usual folding argument to
get $C$ concentrated in two adjacent degrees, in which case the
boundary homomorphism is an $n\delta$ isomorphism.

Recall that we would be finished if we could show that given $\epsilon$
there is $\delta>0$ so that $C$ is $\epsilon$ equivalent to a complex
located over $W^{\epsilon}\cap(B-W^{-\epsilon})$.  The problem is
therefore reduced to a problem about $\delta$ isomorphisms.  Namely we
have a homomorphism $d\:C_j\to C_{j-1}$ that is a $\delta$ isomorphism
off $W^{\delta}\cap(W^{-\delta})$ and we want to split off and discard
a contractible summand containing the part outside
$W^{\epsilon}\cap(B-W^{-\epsilon})$.

It is at this point that our argument diverges from Yamasaki's. 
Splitting an isomorphism is a controlled $K$-theory problem and
generally not possible.  Yamasaki stabilizes by multiplying by $S^1$. 
This canonically kills $K$-theory so splitting becomes possible. 
However it also changes the surgery problem, so he defines
$L^{-\infty}$ by factoring out these changes and concludes that the
assembly is an isomorphism in this context.  In our special case ($\Z$
coefficients) the controlled $K$-theory vanishes and splitting is
possible without stabilization.

\proclaim{7.  Lemma: splitting isomorphisms} Suppose $B\supset V$ and
$\epsilon>0$ are given.  Then there is $\delta>0$ so that if $d\:A\to
A'$ is a $\delta$ homomorphism of $\Z$-modules over $B$ and is a
$\delta$ isomorphism over $B-V$ then there are $\epsilon$ automorphisms
$H$, $H'$ of $A$, $A'$ respectively so that (1) $H$, $H'$ are the
identity over $V^{\epsilon}$, and (2) $H'dH$ is induced by a bijection
of bases over $B-V^{2\epsilon}$.\endproclaim This follows easily from
\cite{Q1, Theorem 8.4}.

We apply this to the homomorphism $d\:C_j\to C_{j-1}$ in the boundary
complex considered above.  Conclusion (1) enables us to extend $H$,
$H'$ by the identity on the rest of $C$ to get an $\epsilon$
equivalence of quadratic Poincar\'e complexes.  Conclusion (2) shows
that the new complex is the sum of a trivial complex ($d=1$) and one
lying over $B-V^{2\epsilon}$.  Deleting the trivial summand gives a
splitting satisfying the conclusions of Lemma 6.  This proves
Lemma 6, which shows that the assembly
$H_n(B;\L)@>>>L_n(B;\Z,\epsilon,\delta)$ is an isomorphism and so
completes the proof of Theorem~1.

\enddocument